\documentclass[12pt,reqno]{article}
\usepackage[a4paper,margin=25mm]{geometry}
\usepackage[T1]{fontenc}
\usepackage{lmodern,microtype}
\usepackage{amsmath,amssymb,amsthm,mathtools}
\usepackage{cite,enumerate}
\usepackage{needspace}
\usepackage[hidelinks]{hyperref}
\newtheorem{theorem}{Theorem}[section]
\newtheorem{lemma}[theorem]{Lemma}
\newtheorem{proposition}[theorem]{Proposition}
 
\newtheorem{remark}{Remark}[theorem]
\numberwithin{equation}{section}

\DeclareMathOperator{\sgn}{sgn}
\DeclareMathOperator{\Impart}{Im}
\newcommand{\RZ}{\mathrm{RZ}}

\newcommand{\RR}{\mathbb R}
\newcommand{\HH}{\mathbb H}

\newcommand{\LL}{\mathcal R}

\newcommand{\dd}{\mathrm d}

\title{Real zeros of $w$-polynomials: cycles and complete cutsets}

\author{Hongyue Tang,
Lily Li Liu
\footnote{{Corresponding author}
\newline\hspace*{5mm}
{\it Email addresses:}
 hytang@qfnu.edu.cn (H. Tang);
 liulily@qfnu.edu.cn (L.L. Liu),}}
\date{\footnotesize
School of Mathematical Sciences,
Qufu Normal University,
Qufu 273165, PR China}
\begin{document}
\maketitle
\begin{abstract}
For a simple graph $G$ with $n$ vertices, its $w$-polynomial is
defined by
\begin{align*}
 \sum_{q\geq0}P(G;q)x^q
 =\frac{w(G;x)}{(1-x)^{n+1}},
\end{align*}
where $P(G;x)$ is the chromatic polynomial of $G$.
This polynomial was introduced by Brenti in 1992.
In this paper, we prove Conjecture~6.7 and the $w$-assertion of
Conjecture~6.5 posed by Brenti, Royle and Wagner in 1994.
The first states that $w(C_n;x)$ has only
real zeros for $n\ge 5$. The second assertion states that if $G\cap H$ is
complete and both $w(G;x)$ and $w(H;x)$ have only real zeros, then
$w(G\cup H;x)$ has only real zeros. 
\end{abstract}
 
\section{Introduction}
\hspace*{\parindent}
Let $G=(V,E)$ be a finite simple graph with $n$ vertices. For a positive integer $q$, a {\it proper $q$-colouring} of $G$ is a map
$c:V(G)\longrightarrow\{1,\ldots,q\}$ such that $c(u)\neq c(v)$
whenever $uv\in E(G)$.
The chromatic polynomial of $G$, denoted by $P(G;x)$, is the
polynomial whose value at each positive integer $q$ counts the
proper $q$-colourings of $G$.
\begin{equation*}
 P(G;x)=\sum_{j=0}^n a_j(G)x_{(j)},
\end{equation*}
where $x_{(j)}=x(x-1)\cdots(x-j+1)$ is the $j$-th falling factorial
polynomial and $a_j(G)$ is the number of partitions of $V(G)$ into
exactly $j$ nonempty independent sets (see \cite{BR1994,BL1946} for instance).

The $w$-polynomial, introduced by Brenti~\cite{B1992} in his study
of chromatic polynomial expansions, is defined by
\begin{equation}\label{E1.1}
 \sum_{q\geq0}P(G;q)x^q
 =\frac{w(G;x)}{(1-x)^{n+1}}.
\end{equation}
 For the graph with no vertices, we use $P(G;x)=w(G;x)=1$.
 
Brenti, Royle and Wagner~\cite{BR1994} subsequently investigated the
zeros of these polynomials and posed Conjectures~6.5 and~6.7.
The following results prove their assertions concerning the $w$-polynomial,
which are the main objective of this paper.
Let $\RZ$ denote the set of real polynomials with only real zeros.

\begin{theorem}[{\rm \cite[Conjecture 6.7]{BR1994}}] \label{T1.1}
Let $C_n$ be a cycle with $n$ vertices. For $n\geq5$, $w(C_n;x) \in \RZ$. Furthermore, the following statements hold.
\begin{enumerate}[(i)]
\item $w(C_5;x)=30x^3(1+x)^2$.
\item If $n\geq6$ is even, zero has multiplicity two and the other
$n-2$ zeros are simple and negative.
\item If $n\geq7$ is odd, zero has multiplicity three and the other
$n-3$ zeros are simple and negative.
\end{enumerate}
\end{theorem}

The following result concerns complete cutsets.

\begin{theorem}[{\rm \cite[Conjecture 6.5]{BR1994}}] \label{T1.2}
Let $G$ and $H$ be two graphs satisfying the following conditions.
\begin{enumerate}[(a)]
\item $G\cap H=K_m$ for some integer $m\geq0$.
\item $w(G;x),w(H;x)\in\RZ$.
\end{enumerate}
Then $w(G\cup H;x)\in\RZ$.
\end{theorem}

In order to prove Theorem~\ref{T1.1}, we first transform the Eulerian recurrence into a differential
recurrence for $R_n(t)$.
Using a derivative inequality obtained by a trigonometric
zero-counting argument, we show that $R_n(t)-1$ and $R_n(t)+1$
have only simple zeros in $[-1,1]$ for $n\geq6$.
The relation between $R_n(t)$ and $A_n(x)$ then gives the
location and multiplicities of the zeros of $w(C_n;x)$. We prove Theorem~\ref{T1.2} by reducing it to an interval-preserving property
of $\diamond_m$. The reduction uses the complete cutset identity for
$\sigma$-polynomials and the operator $\mathcal L_m$. We establish this property by verifying the stability of a
signed symbol using Laguerre polynomials and polarization,
and then excluding zeros outside $[-1,0]$ by a sign analysis.

The organization of this paper is as follows.
In Section~\ref{S2}, we recall the Eulerian identities and introduce
the polynomials $R_n(t)$.
In Section~\ref{S3}, we introduce some properties of trigonometric polynomials and prove Theorem~\ref{T1.1}.
In Section~\ref{S4}, we recall some facts of $\sigma$-polynomials and derive the differential product formula from
the known graph polynomial identities.
In Section~\ref{S5}, we establish the zero-location property of
$\diamond_m$ and prove Theorem~\ref{T1.2}.

\section{Eulerian polynomials and their transformations}\label{S2}
\hspace*{\parindent}
In this section, we recall some facts of the Eulerian polynomials $A_n(x)$ and introduce the operator $\LL$ and a family of polynomials $R_n(t)$. Moreover, we establish the relation between $R_n(t)$ and $A_n(x)$.

Let $\mathfrak{S}_n$ denote the set of permutations on $[n] = \{1,2,\cdots,n \} $.  We say that $i$ is a descent of $\pi = \pi_1 \pi_2 \cdots \pi_n \in \mathfrak{S}_n$ if $1\leq i < n$ and $\pi_i > \pi_{i+1}$. Let des($\pi$) denote the number of descents of $\pi$. The Eulerian polynomials \cite{P2015}, which enumerate the number of descents of $\mathfrak{S}_n$, are defined by

\begin{align}\label{E2.1}
A_n(x) = \sum_{\pi \in \mathfrak{S}_n} x^{\mathrm{des(\pi)}},\quad  A_1(x) = 1.
\end{align}
The Eulerian polynomials and their properties play a pivotal role in combinatorics (see \cite{C2003,LW2007,SV2015} for instance). It is well known that the Eulerian polynomials satisfy
the recurrence relation \cite{SV2012}
\begin{equation}\label{E2.2}
 A_{n+1}(x)=(1+nx)A_n(x)+x(1-x)A_n'(x), \quad n \geq  1,
\end{equation}
and the power-series identity \cite[Proposition~1.4.4]{S2012}
\begin{equation}\label{E2.3}
 \sum_{q\geq0}q^n x^q=\frac{xA_n(x)}{(1-x)^{n+1}}, \quad n \geq  1.
\end{equation}
In order to transform the recurrence into a form that is convenient for
the zero analysis below, we introduce a linear differential operator. For a polynomial $p(t)$,
define
\begin{equation}\label{E2.4}
 \LL p(t)=tp(t)+\frac{t^2-1}{2}p'(t)
   =\frac12\frac{d}{dt}\bigl((t^2-1)p(t)\bigr).
\end{equation}
Since the coefficient of the derivative term
vanishes at $t=\pm1$, the endpoint values can be determined directly
from the recurrence.
Let $\{R_n(t)\}_{n\geq1}$ be the polynomials defined by
\begin{equation}\label{E2.5}
 R_1(t)=1,\qquad R_{n+1}(t)=\LL R_n(t),\qquad n\geq1.
\end{equation}
 
A real polynomial is called \emph{standard} if it is identically
zero or its leading coefficient is positive.

\begin{lemma}\label{L2.1}
For $n\geq1$, we have
\begin{equation}\label{E2.6}
 R_n(t)=\left(\frac{t+1}{2}\right)^{n-1}
           A_n\left(\frac{t-1}{t+1}\right),
\end{equation}
where the right-hand side is interpreted as a polynomial at $t=-1$.
Moreover, $R_n(t)$ is standard of degree $n-1$ and
\begin{equation}\label{E2.7}
 R_n(-1)=(-1)^{n-1},\qquad R_n(1)=1.
\end{equation}
\end{lemma}

\begin{proof}
For $t\neq-1$, let
\begin{align*}
 x(t)=\frac{t-1}{t+1},
 \qquad
 B_n(t)=\frac{A_n(x(t))}{(1-x(t))^{n-1}}.
\end{align*}
We first show that $B_n(t)$ is a polynomial in $t$.
Suppose that
\begin{align*}
 A_n(x)=\sum_{j=0}^{n-1}a_{n,j}x^j,
\end{align*}
where $a_{n,j}$ is the number of permutations of $n$ elements
with exactly $j$ descents.

Since $ 1-x(t) =\frac{2}{t+1}$, we obtain
\begin{align*}
 B_n(t)
 &=\left(\frac{t+1}{2}\right)^{n-1}
   \sum_{j=0}^{n-1}a_{n,j}
       \left(\frac{t-1}{t+1}\right)^j\\
 &=2^{-(n-1)}
   \sum_{j=0}^{n-1}
       a_{n,j}(t-1)^j(t+1)^{n-1-j}.
\end{align*}
Hence $B_n(t)$ is a polynomial
and its value at $t=-1$ is also well defined.

We next derive a recurrence for $B_n(t)$.
For $t\neq-1$,  by
differentiation we obtain
\begin{align*}
 B_n'(t)
 &=\frac{x'(t)A_n'(x(t))}{(1-x(t))^{n-1}}
   +\frac{(n-1)x'(t)A_n(x(t))}{(1-x(t))^n}\\
 &=x'(t)\left(
       \frac{A_n'(x(t))}{(1-x(t))^{n-1}}
       +\frac{n-1}{1-x(t)}B_n(t)\right).
\end{align*}
 Hence, we have
\begin{align} \label{E2.8}
 \frac{A_n'(x(t))}{(1-x(t))^{n-1}}
 =\frac{B_n'(t)}{x'(t)}
  -\frac{n-1}{1-x(t)}B_n(t).
  \end{align}
By equalities \eqref{E2.2} and \eqref{E2.8}, we conclude that
\begin{align*}
 B_{n+1}(t)
 &=\frac{A_{n+1}(x(t))}{(1-x(t))^n}\\
 &=\frac{(1+nx(t))A_n(x(t))+x(t)(1-x(t))A_n'(x(t))}
         {(1-x(t))^n}\\
 &=\frac{1+nx(t)}{1-x(t)}B_n(t)
   +x(t)\left(
       \frac{B_n'(t)}{x'(t)}
       -\frac{n-1}{1-x(t)}B_n(t)\right)\\
 &=\frac{1+x(t)}{1-x(t)}B_n(t)
   +\frac{x(t)}{x'(t)}B_n'(t).
\end{align*}
By direct calculation, we have
\begin{align*}
 \frac{1+x(t)}{1-x(t)} =t \quad \mathrm{and} \quad
 \frac{x(t)}{x'(t)} =\frac{t^2-1}{2}.
\end{align*}
Hence, 
\begin{align*}
 B_{n+1}(t)
 =tB_n(t)+\frac{t^2-1}{2}B_n'(t)
 =\LL B_n(t),
 \qquad t \neq-1.
\end{align*}
Both sides are polynomials in $t$, so their difference vanishes at
all real $t\neq-1$, and hence has infinitely many zeros.
It must be the zero polynomial. Thus the
recurrence also holds at $t=-1$.

We next show that $R_n(t)=B_n(t)$.
Since $A_1(x)=1$, the polynomial $ B_1(t)=1=R_1(t)$.
Suppose that $B_n(t)=R_n(t)$ for some $n\geq1$.
Combining the recurrence above with equality \eqref{E2.5}, we obtain
 \begin{align*}
 B_{n+1}(t)
 =\LL B_n(t)
 =\LL R_n(t)
 =R_{n+1}(t).
\end{align*}
It follows by induction that $B_n(t)=R_n(t)$ for $n\geq1$.
Thus we have
\begin{align*}
 R_n(t)
 =\frac{A_n(x(t))}{(1-x(t))^{n-1}}
 =\left(\frac{t+1}{2}\right)^{n-1}
   A_n\left(\frac{t-1}{t+1}\right),
\end{align*}
where the right-hand side is a polynomial in $t$
and hence is defined at $t=-1$.

It remains to determine the degree, the sign of the leading coefficient,
and the endpoint values of $R_n(t)$.
The polynomial $R_1(t)=1$ has degree zero and positive
leading coefficient.
Suppose that $R_n(t)$ has degree $n-1$ and leading coefficient
$c_n>0$. By the recurrence
\begin{align*}
 R_{n+1}(t)
 =tR_n(t)+\frac{t^2}{2}R_n'(t)-\frac12R_n'(t),
\end{align*}
the term $tR_n(t)$ contributes $c_n$ to the coefficient of $t^n$,
and the term $\frac{t^2}{2}R_n'(t)$ contributes $(n-1)c_n/2$.
The term $-\frac12R_n'(t)$ contributes zero.
Thus the coefficient of $t^n$ in $R_{n+1}(t)$ is
\begin{align*}
 c_n+\frac{n-1}{2}c_n
 =\frac{n+1}{2}c_n>0.
\end{align*}
Hence $R_{n+1}(t)$ has degree $n$ and positive leading coefficient.
By induction, $R_n(t)$ is standard of degree $n-1$.

Finally, by substituting $t=1$ and $t=-1$ into the recurrence, we obtain
\begin{align*}
 R_{n+1}(1) =R_n(1)=1 \quad \mathrm{and} \quad
 R_{n+1}(-1)=-R_n(-1)= (-1)^n.
\end{align*}
 
This completes the proof.

\end{proof}

\section{The real zeros of \texorpdfstring{$w(C_n;x)$}{w(Cn;x)}}\label{S3}
\hspace*{\parindent}
In this section we prove that $w(C_n;x)$ has only real zeros
for $n\geq5$. We first introduce trigonometric polynomials and show that $R_n(t)-1$ and $R_n(t)+1$ have only simple zeros
in $[-1,1]$ for $n\geq6$. By a change of variables, we obtain the factorization of $w(C_n;x)$
and the zero multiplicities stated in Theorem~\ref{T1.1}.

A real trigonometric polynomial of degree at most $d$ is an expression of the form
\begin{equation*}
 Q(\theta)=a_0+\sum_{k=1}^d
 \bigl(a_k\cos(k\theta)+b_k\sin(k\theta)\bigr),
\end{equation*}
where $a_k,b_k\in\mathbb{R}$.
The following result provides a bound on the number of its zeros.
\begin{lemma}\label{L3.1}
Let $d\geq1$ and let $Q(\theta)$ be a nonzero real trigonometric polynomial of
degree at most $d$. Then $Q(\theta)$ has at most $2d$ distinct zeros in a
half-open interval of length $2\pi$.
\end{lemma}
\begin{proof}
Let $z=e^{i\theta}=\cos\theta+i\sin\theta$. Using the identities
$\cos(k\theta)=(z^k+z^{-k})/2$ and
$\sin(k\theta)=(z^k-z^{-k})/(2i)$, we obtain
\begin{align*}
 z^dQ(\theta)=a_0z^d+
 \sum_{k=1}^d\left(
     \frac{a_k-ib_k}{2}z^{d+k}
    +\frac{a_k+ib_k}{2}z^{d-k}\right).
\end{align*}
It is clear that $z^dQ(\theta)$ is a nonzero ordinary polynomial of degree at most $2d$. If $Q(\theta)=0$, its associated nonzero number
$z=e^{i\theta}$ is a zero of this polynomial.
Since distinct angles in a half-open interval of length $2\pi$ give distinct nonzero values of $z$,
the result follows from the bound on the number of zeros of an ordinary polynomial.
\end{proof}

The following lemma provides a lower bound on the derivative of a
trigonometric polynomial at points where its value lies strictly
between $-1$ and $1$.
\begin{lemma}\label{L3.2}
Let $Q(\theta)$ be a real trigonometric polynomial of degree at most $d\geq1$.
Suppose that the following conditions are satisfied.
\begin{enumerate}[(a)]
\item There are points
$\theta_0<\theta_1<\cdots<\theta_{2d}=\theta_0+2\pi$ such that
$(-1)^jQ(\theta_j)\geq1$ for $0\leq j\leq2d$.
\item $Q(\theta)$ is strictly monotone on each
$(\theta_j,\theta_{j+1})$.
\end{enumerate}
If $|Q(\varphi)|<1$ and $Q'(\varphi)\ne0$, then
\begin{equation}\label{E3.1}
 Q'(\varphi)^2\geq d^2\big(1-Q(\varphi)^2\big).
\end{equation}
\end{lemma}
\begin{proof}
Suppose, to the contrary, that
\begin{equation*}
 |Q'(\varphi)|<d\sqrt{1-Q(\varphi)^2}.
\end{equation*}
By periodicity,
we may assume that $\varphi\in(\theta_j,\theta_{j+1})$ for some $j\in\{0,1,\ldots,2d-1\}$. Let
$u=Q(\varphi)$ and $w=Q'(\varphi)$.
Choose a real number $v$ with the same sign as $Q'(\varphi)$ such that
\begin{align}\label{E3.2}
|w|<|v|<d\sqrt{1-u^2}.
\end{align}
Define
\begin{align*}
 C(\theta) =u\cos(d(\theta-\varphi))
             +\frac vd\sin(d(\theta-\varphi))\quad \mathrm{and} \quad S(\theta)=Q(\theta)-C(\theta).
\end{align*}
By direct calculation, we have
\begin{align*}
 C(\varphi)=u= Q(\varphi),\qquad C'(\varphi)=v,\qquad
 S(\varphi)=0,\qquad S'(\varphi)=w-v.
\end{align*}
By inequality \eqref{E3.2} and the Cauchy--Schwarz inequality, we obtain

\begin{align*}
 |C(\theta)|
 & \leq \sqrt{u^2+\frac{v^2}{d^2}} \cdot
       \sqrt{\cos^2(d(\theta-\varphi))+\sin^2(d(\theta-\varphi))}\\
 &= \sqrt{u^2+\frac{v^2}{d^2}} <1.
\end{align*}
Hence $(-1)^jC(\theta_j) \in (-1,1)$ and
\begin{align*}
 (-1)^jS(\theta_j)
 =(-1)^jQ(\theta_j)-(-1)^jC(\theta_j)>0.
\end{align*}
In particular, $S(\theta)$ is not identically zero and $S(\theta)$ changes sign in
each of the $2d$ open intervals $(\theta_j,\theta_{j+1})$.

By periodicity, let $\varphi \in [\theta_0,\theta_0+2\pi)$.
It cannot be one of the $\theta_j$, since $|Q(\varphi)|<1$.
Assume that $\varphi\in(\theta_j,\theta_{j+1})$ and $\sigma=(-1)^j$.  
At the left endpoint $\sigma Q(\theta_j)= (-1)^{j}Q(\theta_j) \geq1$, and at the right endpoint
$\sigma Q(\theta_{j+1})  =(-1)^{j} Q(\theta_{j+1}) \leq-1$. By monotonicity, $\sigma Q(\theta)$ is strictly decreasing.
Since $w\ne0$, it follows that $\sigma w<0$.
The choice of $v$ then gives $\sigma v<\sigma w<0$ and hence
\begin{align}\label{E3.3}
 \sigma S'(\varphi)=\sigma(w-v)>0.
\end{align}
By inequality \eqref{E3.3} and the equality $S(\varphi)=0$, there is an $\varepsilon>0$ such that
\begin{align*}
 \sigma S(\varphi-\varepsilon)<0 \quad {\rm and} \quad
 \sigma S(\varphi+\varepsilon)>0.
\end{align*}
On the other hand,
$\sigma S(\theta_j)>0$ and $\sigma S(\theta_{j+1})<0$.
There is a zero between
$\theta_j$ and $\varphi-\varepsilon$, and another zero between
$\varphi+\varepsilon$ and $\theta_{j+1}$ by the intermediate value theorem.
Thus, these are three distinct zeros in the same
open interval. Each of the remaining $2d-1$ open intervals has at
least one zero. All these intervals are disjoint, so the number
of distinct zeros in one period is at least
$3+(2d-1)=2d+2$.
By Lemma~\ref{L3.1}, $S(\theta)$ has at most $2d$ zeros. This
contradicts the result that $S$ has $2d+2$ zeros. Hence
inequality~\eqref{E3.1} holds.
\end{proof}

Since Lemma~\ref{L3.2} requires monotonicity intervals of
$Q(\theta)=p(\cos\theta)$ whose endpoint values alternate beyond
$\pm1$, we prove the following lemma to construct such intervals.

\begin{lemma}\label{L3.3}
Let $p(t)$ be a standard polynomial of degree $d\geq1$ satisfying:
\begin{enumerate}[(a)]
\item $p(-1)=(-1)^d$ and $p(1)=1$;
\item $p(t)-1$ and $p(t)+1$ have $d$ distinct zeros in $[-1,1]$, respectively.
\end{enumerate}
Then the following results hold:
\begin{enumerate}[(i)]
\item If $d\geq2$, each local maximum of $p(t)$ is greater than $1$
and each local minimum is less than $-1$.

\item For each $s\in(-1,1)$, the polynomial $p(t)-s$ has $d$ distinct
zeros in $(-1,1)$.

\item There are ordered disjoint open intervals $I_1,I_2,\ldots,I_d$
in $(-1,1)$ on which $p(t)$ maps strictly monotonically onto $(-1,1)$,
and the sign of $p'(t)$ on $I_i$ is $\sigma_i=(-1)^{d-i}$.
\end{enumerate}
\end{lemma}

\begin{proof}
If $d=1$, we have $p(t)=t$ by the endpoint condition $(a)$ and the result
follows immediately.
Assume that $d\geq2$ and let $a>0$ be the leading coefficient of $p(t)$. Let
\begin{align*}
 u_1<u_2<\cdots<u_d\quad {\rm and} \quad v_1<v_2<\cdots<v_d
\end{align*}
be the zeros of $p(t)-1$ and $p(t)+1$, respectively.
By Rolle's theorem, $p'(t)$ has a zero in each interval
$(u_i,u_{i+1})$ where $1\leq i\leq d-1$. Hence $p'(t)$ has $d-1$ distinct zeros in $(-1,1)$.
Since its degree is $d-1$, these are all its zeros and each is simple.
Suppose that the zeros of $p'(t)$ are $c_1<c_2<\cdots<c_{d-1}$.
The result also holds for $p(t)+1$. So we obtain
 \begin{align}\label{E3.4}
 u_i<c_i<u_{i+1} \quad {\rm and} \quad
  v_i<c_i<v_{i+1},\qquad(1\leq i\leq d-1).
\end{align}

Under the hypothesis, $ p(t)-1$ and $ p(t)+1$ can be represented in the following forms, respectively.
\begin{align*}
 p(t)-1=a\prod_{j=1}^d(t-u_j),\qquad
 p(t)+1=a\prod_{j=1}^d(t-v_j).
\end{align*}
By inequality \eqref{E3.4}, $p(c_i)-1$ and $p(c_i)+1$ have exactly $d-i$ negative factors, respectively.
Their signs are given by
\begin{align}\label{E3.5}
 \operatorname{sgn}\bigl(p(c_i)-1\bigr)
 =\operatorname{sgn}\bigl(p(c_i)+1\bigr)
 =(-1)^{d-i}.
\end{align}
Let $c_0=-1$ and $c_d=1$. Since
\begin{align*}
 p'(t)=da\prod_{j=1}^{d-1}(t-c_j),
\end{align*}
there are exactly $d-i$ negative factors when $c_{i-1}<t<c_i$. Hence
\begin{align*}
 \operatorname{sgn}p'(t)=(-1)^{d-i},
 \qquad t\in(c_{i-1},c_i),\quad 1\leq i\leq d.
\end{align*}
The derivative $p'(t)$ changes from positive to negative at $t=c_i$
if $(-1)^{d-i}=1$, and from negative to positive if $(-1)^{d-i}=-1$.
Hence, by equality~\eqref{E3.5}, we have $p(c_i)>1$ for the case of $(-1)^{d-i}=1$ and $p(c_i)<-1$ for the case of $(-1)^{d-i}=-1$.
Thus each local maximum of $p(t)$ is greater than $1$ and each local minimum
is less than $-1$.

For any $s\in(-1,1)$, by equality~\eqref{E3.5} and condition $(a)$, we have 
\begin{align*}
 \operatorname{sgn}\bigl(p(c_i)-s\bigr)=(-1)^{d-i},
 \qquad 0\leq i\leq d.
\end{align*}
By the intermediate value theorem, it has a zero $r_i(s)$
in each interval $(c_{i-1},c_i)$ where $1 \leq i \leq d$.
These are $d$ distinct zeros in $(-1,1)$.
Since $\deg(p(t)-s)=d$, we have
\begin{align*}
 p(t)-s=a\prod_{i=1}^d\bigl(t-r_i(s)\bigr).
\end{align*}
Thus these are all its zeros and each is simple.

It remains to construct the intervals $I_i$.
For $1 \leq i \leq d$, let $J_i=(c_{i-1},c_i)$ and $\sigma_i=(-1)^{d-i}$.
The derivative sign obtained above shows that $\sigma_i p(t)$
is strictly increasing on $[c_{i-1},c_i]$.
In particular, we have
\begin{align*}
 \sigma_i p(c_{i-1})\leq-1,\qquad
 \sigma_i p(c_i)\geq1.
\end{align*}
Hence, by continuity and strict monotonicity, there are unique points $\alpha_i$ and $\beta_i$  satisfying
\begin{align*}
 c_{i-1}\leq\alpha_i<\beta_i\leq c_i,\qquad
 \sigma_i p(\alpha_i)=-1,\qquad
 \sigma_i p(\beta_i)=1.
\end{align*}
Define $I_i=(\alpha_i,\beta_i)$.
The boundary values and strict monotonicity imply that
$\sigma_i p(t)$ maps $I_i$ onto $(-1,1)$.
Since multiplication by $\sigma_i\in\{-1,1\}$ preserves this
range, $p(t)$ also maps $I_i$ strictly monotonically onto $(-1,1)$.

Since $I_i\subseteq J_i$, the intervals $I_1,\ldots,I_d$
are ordered, disjoint and contained in $(-1,1)$.
The sign of $p'(t)$ on $I_i$ is
\begin{align*}
 \operatorname{sgn}p'(t)=\sigma_i=(-1)^{d-i},
 \qquad t\in I_i.
\end{align*}
This completes the proof.

\end{proof}
 
 The critical points and monotonicity intervals supplied above allow us
to apply the trigonometric derivative bound to \(Q(\theta)=p(\cos\theta)\).
Hence, we obtain the following result.

\begin{lemma}\label{L3.4}
Suppose that $p(t)$ satisfies the conditions of Lemma~\ref{L3.3}.
For any $u\in(-1,1)$ such that $|p(u)|<1$, we have
\begin{equation}\label{E3.6}
 (1-u^2)p'(u)^2\geq d^2\bigl(1-p(u)^2\bigr).
\end{equation}
\end{lemma}
\begin{proof}
For $d=1$, $p(t)=t$, so both sides equal $1-u^2$.
Assume that $d\geq2$, and let $Q(\theta)=p(\cos\theta)$.
For integers $k$ with $0\leq k\leq d$, we have
\begin{align*}
 (\cos\theta)^k
 =2^{-k}(e^{i\theta}+e^{-i\theta})^k
 =2^{-k}\sum_{\ell=0}^k\binom{k}{\ell}
                 e^{i(k-2\ell)\theta}.
\end{align*}
Thus $Q(\theta)$ is a real trigonometric polynomial of degree at most $d$. Next, we construct an increasing sequence of points. Let $(c_j)_{j=1}^{d-1}$ be the zeros of $p'(t)$ as in the proof of Lemma~\ref{L3.3}.
 For $1\leq j<d $, define
\begin{align*}
\begin{split}
 \theta_0=0,\ldots, \theta_j=\arccos(c_{d-j}),\ldots, \theta_d=\pi,\ldots, \theta_{d+j}=2\pi-\arccos(c_j),\ldots, \theta_{2d}=2\pi.
\end{split}
\end{align*}
Since $c_1<c_2<\cdots<c_{d-1}$ and $\arccos$ is strictly decreasing,
the points $(\theta_i)_{i=0}^{2d}$ form an increasing
sequence.
By equality~\eqref{E3.5}, we have
\begin{align*}
 \sgn Q(\theta_j) = \sgn p(c_{d-j}) =(-1)^j,\qquad |Q(\theta_j)|>1,
\end{align*}
and 
\begin{align*}
 \sgn Q(\theta_{d+j}) = \sgn p(c_j) =(-1)^{d-j}= (-1)^{d-j} (-1)^{-2j}=(-1)^{d+j},\qquad |Q(\theta_{d+j})|>1.
\end{align*}
The remaining points are $Q(0)=Q(2\pi)=1$ and $Q(\pi)=(-1)^d$.
Hence $(-1)^iQ(\theta_i)\geq1$ for $0 \leq i \leq 2d$.
Since $ Q'(\theta)=-\sin\theta\,p'(\cos\theta)$ and between consecutive $\theta_j$ neither $\sin\theta$ nor
$p'(\cos\theta)$ vanishes, $Q(\theta)$ is strictly monotone on each interval $(\theta_i,\theta_{i+1})$.

Let $\varphi=\arccos u\in(0,\pi)$.
Since $|p(u)|<1$,
$\varphi$ lies strictly between consecutive cyclic points and $Q'(\varphi)\neq0$.
By Lemma~\ref{L3.2} and $\sin^2\varphi=1-u^2$,
we obtain
\begin{align*}
 (1-u^2)p'(u)^2
 =Q'(\varphi)^2
 \geq d^2(1-Q(\varphi)^2)
 =d^2(1-p(u)^2).
\end{align*}
This completes the proof.
 
\end{proof}

We now return to the operator $\LL$. By Lemmas~\ref{L3.3}
and~\ref{L3.4}, we control $\LL p(t)$ at
suitable points and prove that the polynomials $R_n(t)$ constructed by
iterating $\LL$ preserve the property that $R_n(t)-1$ and
$R_n(t)+1$ have only simple zeros in $[-1,1]$. 

\begin{theorem}\label{T3.5}
Let $p(t)$ be a standard polynomial of degree $d\geq3$ satisfying:
\begin{enumerate}[(a)]
\item $p(-1)=(-1)^d$ and $p(1)=1$;
\item $p(t)-1$ and $p(t)+1$ have $d$ distinct zeros in $[-1,1]$, respectively.
\end{enumerate}
Then the following results hold:
\begin{enumerate}[(i)]
\item $\LL p(t)$ is a standard polynomial of degree $d+1$.

\item $\LL p(-1) =(-1)^{d+1}$ and $\LL p(1)=1$.

\item $\LL p(t)-1$ and $\LL p(t)+1$ have $d+1$ distinct zeros in $[-1,1]$, respectively.
\end{enumerate}
\end{theorem}
\begin{proof}
The degree, the leading coefficient and the endpoint values follow directly from equality \eqref{E2.4}.
So it suffices to show that $\LL p(t)-1$ and $\LL p(t)+1$ have $d+1$ distinct zeros in $[-1,1]$, respectively.

By Lemma~\ref{L3.3},
there are ordered intervals $I_i=(\alpha_i,\beta_i)$,
where $-1\leq\alpha_i<\beta_i\leq1$,
such that $p(t)$ maps each $I_i$ strictly monotonically onto $(-1,1)$.
Let $\sigma_i=(-1)^{d-i}$ be the sign of $p'(t)$ on $I_i$.

(i) If $\sigma_i=1$, then $p(\alpha_i)=-1$ and $p(\beta_i)=1$.
Hence
\begin{equation*}
 p(\alpha_i)+\alpha_i=-1+\alpha_i<0
 \quad\mathrm{and}\quad
 p(\beta_i)+\beta_i=1+\beta_i>0.
\end{equation*}
By the intermediate value theorem,
there is a point $u_i\in I_i$ such that $p(u_i)=-u_i$.

(ii) If $\sigma_i=-1$, then $p(\alpha_i)=1$ and $p(\beta_i)=-1$.
Thus
\begin{equation*}
 p(\alpha_i)-\alpha_i=1-\alpha_i>0
 \quad\mathrm{and}\quad
 p(\beta_i)-\beta_i=-1-\beta_i<0.
\end{equation*}
By the intermediate value theorem, there is a point $u_i\in I_i$ such that $p(u_i)=u_i$.

In both cases,
we have
\begin{equation}\label{E3.7}
 -1<u_1<u_2<\cdots<u_d<1
 \quad\mathrm{and}\quad
 p(u_i)=-\sigma_i u_i.
\end{equation}
The endpoint inequalities above are strict since $\alpha_i<1$ and $\beta_i>-1$. By Lemma~\ref{L3.4} and \eqref{E3.7}, we obtain
\begin{align*}
 (1-u_i^2)p'(u_i)^2
 \geq d^2(1-p(u_i)^2)=d^2(1-u_i^2).
\end{align*}
Since $|u_i|<1$,
it follows that $|p'(u_i)|\geq d$.
Moreover, by Lemma \ref{L3.3} (iii), $p'(u_i)=\sigma_i|p'(u_i)|$.
Thus
\begin{align*}
 \LL p(u_i)
 &=u_i(-\sigma_i u_i)
       +\frac{u_i^2-1}{2}\sigma_i|p'(u_i)|\\
 &=-\sigma_i\left(u_i^2+
                \frac{1-u_i^2}{2}|p'(u_i)|\right).
\end{align*}
Multiplying by $-\sigma_i$ and applying $|p'(u_i)|\geq d$, we have
\begin{align}
 -\sigma_i \LL p(u_i)
 &\geq u_i^2+\frac d2(1-u_i^2)\notag\\
 &=1+\left(\frac d2-1\right)(1-u_i^2)>1.\label{E3.8}
\end{align}
The last inequality uses $d\geq3$ and $|u_i|<1$. Thus the values $\LL p(u_i)$ have alternating signs and absolute value greater than $1$.

For any $s\in\{-1,1\}$,
the values of $\LL p(t) -s$ at consecutive points $u_i$ and $u_{i+1}$ have opposite signs.
So there are $d-1$ distinct zeros between these points.
There is one further zero in $[-1,u_1)$. If $\LL p(-1)=s$, then $-1$ is a zero.
If $\LL p(-1) \neq s$, then $\LL p(-1)-s$ and $\LL p(u_1)-s$ have opposite signs.
Similarly,
there is one further zero in $(u_d,1]$.
If $s=1$, then $1$ is a zero.
If $s=-1$, then $\LL p(u_d)+1<0< \LL p(1)+1$.
Hence $\LL p(t)-s$ has $d+1$ distinct zeros in $[-1,1]$.
Since its degree is $d+1$,
these are all its zeros and all of them are simple.
This completes the proof.
\end{proof}

\begin{lemma}\label{L3.6}
For $n\geq6$, $R_n-1$ and $R_n+1$ have $n-1$ distinct zeros
in $[-1,1]$.
\end{lemma}
\begin{proof}
By direct calculations from equality~\eqref{E2.5},
we obtain
\begin{equation*}
 R_6(t)=\frac{1}{2}(45t^5-60t^3+17t).
\end{equation*}
Its values at the following six rational points are
\begin{equation*}
\begin{array}{c|rrrrrr}
 t&-1&-4/5&-1/3&1/3&4/5&1\\ \hline
R_6(t)&-1&742/625&-49/27&49/27&-742/625&1.
\end{array}
\end{equation*}
Thus $R_6(t)-1$ has a zero in each of
\begin{equation*}
 (-1,-4/5),\quad(-4/5,-1/3),\quad
 (-1/3,1/3),\quad(1/3,4/5),
\end{equation*}
and it has a further zero at $t=1$.
Similarly,
$R_6(t)+1$ has a zero in each of
\begin{equation*}
 (-4/5,-1/3),\quad(-1/3,1/3),\quad
 (1/3,4/5),\quad(4/5,1),
\end{equation*}
and it has a further zero at $t=-1$.
In each case these are five distinct zeros in $[-1,1]$.
Since the degree is five,
these are all the zeros and all of them are simple.
The endpoint factors are given by
\begin{align*}
 R_6(t)-1&=\frac{t-1}{2}
      (45t^4+45t^3-15t^2-15t+2),\\
 R_6(t)+1&=\frac{t+1}{2}
      (45t^4-45t^3-15t^2+15t+2).
\end{align*}

Assume that the result holds for some $n\geq6$.
By equality~\eqref{E2.7},
$R_n(t)$ satisfies the endpoint conditions in Theorem~\ref{T3.5}.
Since $\deg R_n(t)=n-1\geq5$ and $R_{n+1}(t)=\LL R_n(t)$, $R_{n+1}(t)-1$ and $R_{n+1}(t)+1$ each have $n$ distinct zeros in $[-1,1]$ by Theorem~\ref{T3.5}.
The result follows by induction. This completes the proof.
\end{proof}

Using the properties of the operator $\LL$ derived above, we prove one of the main results of this paper as follows.
\begin{proof}[Proof of Theorem~\ref{T1.1}]
The chromatic polynomial of a cycle $C_n$ with $n$ vertices is
\begin{equation}\label{E3.9}
 P(C_n;k)=(k-1)^n+(-1)^n(k-1).
\end{equation}

For the first term in equality~\eqref{E3.9}, separating the term $k=0$
and replacing $k-1$ by $m$ in the remaining sum, we have
\begin{align*}
 \sum_{k\geq0}(k-1)^nx^k
 =(-1)^n+x\sum_{m \geq0}m^nx^m
 =(-1)^n+\frac{x^2A_n(x)}{(1-x)^{n+1}}.
\end{align*}
For the second term, the geometric series and its derivative give
\begin{align*}
 \sum_{k\geq0}(k-1)x^k
 =\frac{x}{(1-x)^2}-\frac1{1-x}
 =\frac{2x-1}{(1-x)^2}.
\end{align*}
Combining these identities with~\eqref{E1.1},
we have
\begin{align*}
 w(C_n;x)
 &=x^2A_n(x)+(-1)^n(1-x)^{n+1}
       +(-1)^n(2x-1)(1-x)^{n-1}\\
 &=x^2A_n(x)+(-1)^n(1-x)^{n-1}
                         \big((1-x)^2+2x-1\big)\\
 &=x^2\big(A_n(x)+(-1)^n(1-x)^{n-1}\big).
\end{align*}
Suppose that
\begin{equation}\label{E3.10}
 Q_n(x)=A_n(x)+(-1)^n(1-x)^{n-1} \quad {\rm and}
 \quad w(C_n;x)=x^2Q_n(x).
\end{equation}

For $n=5$,
by direct calculations we obtain $Q_5(x)=30x(1+x)^2$.
Thus (i) holds.

For $n\geq6$, let  
\begin{equation*}
H_n(t)=R_n(t)+(-1)^n.
\end{equation*}
By Lemma~\ref{L3.6},
$H_n(t)$ has $n-1$ distinct zeros $r_1 <r_2<\cdots<r_{n-1}$ in $[-1,1]$. Set $t=\frac{1+x}{1-x}$. By equality \eqref{E2.6}, we have 
\begin{align*}
A_n(x)=(1-x)^{n-1}R_n \bigg(\frac{1+x}{1-x}\bigg)
\end{align*}
and
\begin{align}
 Q_n(x) &=(1-x)^{n-1}R_n \bigg(\frac{1+x}{1-x}\bigg)+(-1)^n(1-x)^{n-1} \notag \\
 &= (1-x)^{n-1}H_n\bigg(\frac{1+x}{1-x}\bigg). \label{E3.11}
\end{align}
This equality is first obtained for $x\neq1$
and is a polynomial identity after clearing the denominators.

Let $a>0$ be the leading coefficient of $H_n(t)$.
By equality~\eqref{E3.11},
we obtain
\begin{align}
 Q_n(x)
 &=(1-x)^{n-1} a\prod_{j=1}^{n-1}
             \left(\frac{1+x}{1-x}-r_j\right)\notag\\
 &=a\prod_{j=1}^{n-1}
             \big(1+x-r_j(1-x)\big)\notag\\
 &=a\prod_{j=1}^{n-1}
             \big((1-r_j)+(1+r_j)x\big).\label{E3.12}
\end{align}
By equality \eqref{E2.7}, $H_n(-1)=(-1)^{n-1}+(-1)^n=0$.
This is a simple zero by Lemma~\ref{L3.6}. The corresponding factor in equality~\eqref{E3.12} is the nonzero constant $2$, so it contributes no finite zero.

For $1<j\leq n-1$, we have $r_j\in(-1,1]$, and the corresponding factor has a simple zero
\begin{equation*}
 x_j=\frac{r_j-1}{r_j+1}.
\end{equation*}
These zeros are distinct because the function $(t-1)/(t+1)$ is strictly
increasing on $(-1,1]$, with derivative $2/(t+1)^2>0$.
 Note that $H_n(1)=1+(-1)^n$.
So $1$ is a zero of $H_n(t)$ if and only if $n$ is odd
and its simplicity follows from Lemma~\ref{L3.6}. By differentiating~\eqref{E2.5},
we obtain
\begin{equation*}
 R_{n+1}'(t)=R_n(t)+2tR_n'(t)
              +\frac{t^2-1}{2}R_n''(t).
\end{equation*}
Thus $R_{n+1}'(1)=1+2R_n'(1)$.
Since $R_1'(1)=0$,
we have
\begin{equation*}
 H_n'(1)=R_n'(1)=2^{n-1}-1>0,
 \quad n\geq2.
\end{equation*}
Since $r_1=-1$ and $r_j>-1$ for $j\geq2$, equality~\eqref{E3.12}
contains exactly one constant factor and $n-2$ nonconstant linear factors.
Thus $\deg Q_n(x)=n-2$. If $n$ is odd,
the factor corresponding to $r_{n-1}=1$ is $2x$, which contributes the zero $x=0$ to $Q_n(x)$.
Every other factor in~\eqref{E3.12} has nonzero constant term. Since $r_1 =-1$ does not contribute a zero to $Q_n(x)$, the remaining $n-3$ zeros $r_j$ of 
$H_n(t)$ yield the remaining $n-3$
distinct zeros $(r_j-1)/(r_j+1)$ of $Q_n(x)$.
 
If $n$ is even,
all $n-2$ zeros of $Q_n(x)$ are simple and negative.
 By equality \eqref{E3.10}, we obtain (ii) and (iii).
This completes the proof.
\end{proof}

\begin{remark}
By direct calculations,
we have
\begin{equation*}
 w(C_4;x)=2x^2(1+4x+7x^2).
\end{equation*}
It has complex zeros $ x=\frac{1}{7}(-2\pm i\sqrt3)$,
so $w(C_4;x)$ does not have only real zeros.
\end{remark}

\section{ Graph polynomials and complete cutsets}\label{S4}
\hspace*{\parindent}
In this section, we recall the $\sigma$-polynomial, the augmented
$\sigma$-polynomial and their relation with the $w$-polynomial.
Then we use the complete cutset identity of Brenti, Royle and
Wagner to obtain a differential product formula.

Recall that $a_k(G)$ denotes the number of partitions of $V(G)$
into exactly $k$ nonempty independent sets.
Following~\cite{BR1994}, the $\sigma$-polynomial and the augmented
$\sigma$-polynomial of $G$ are defined by
\begin{align*}
 \sigma(G;x)=\sum_{k=0}^n a_k(G)x^k
\end{align*}
and
\begin{align*}
 \bar{\sigma}(G;x)=\sum_{k=0}^n k!a_k(G)x^k,
\end{align*}
respectively.
The following lemma describes the relation between the
\(w\)-polynomial and the augmented \(\sigma\)-polynomial.
 
\begin{lemma}[{\cite{BR1994}}]\label{L4.1}
Let $G$ be a graph with $n$ vertices. Then
\begin{align}
 w(G;z)
 &=(1-z)^n\bar{\sigma}\bigg(G;\frac{z}{1-z}\bigg),
 \label{E4.2}\\
 \bar{\sigma}(G;x)
 &=(1+x)^nw\bigg(G;\frac{x}{1+x}\bigg).
 \label{E4.3}
\end{align}
Moreover,
\begin{equation}\label{E4.4}
 w(G;x)\in\RZ
 \quad\Longleftrightarrow\quad
 \text{all zeros of }\bar{\sigma}(G;x)
 \text{ belong to }[-1,0].
\end{equation}
\end{lemma}

We use the $*$-product introduced by Brenti, Royle and
Wagner~\cite{BR1994}. Let $D=\mathrm d/\mathrm dx$. For
$p(x),q(x)\in\RR[x]$, this product is defined by
\begin{equation}\label{E4.5}
 p(x)\ast q(x)=\sum_{k\geq0}\frac{x^k}{k!}(D^kp(x))(D^kq(x)).
\end{equation}
Since $D^kp(x)=D^kq(x)=0$ for all sufficiently large $k$, the sum is
finite.
The following theorem provides the complete cutset identity for
$\sigma$-polynomials.
\begin{theorem}[\cite{BR1994}]\label{T4.2}
Let $G$ and $H$ be graphs such that $G\cap H$ is complete, and let
$m=|V(G\cap H)|$. Then
\begin{equation}\label{E4.6}
 \frac{\sigma(G\cup H;x)}{x^m}
 =
 \frac{\sigma(G;x)}{x^m}
 \ast
 \frac{\sigma(H;x)}{x^m}.
\end{equation}
\end{theorem}

 We next transfer the product $\ast$ to the
augmented $\sigma$-polynomials.
For an integer $m\geq0$, define a linear operator
$\mathcal L_m:\RR[x]\longrightarrow\RR[x]$ by
\begin{equation}\label{E4.7}
 \mathcal L_m(x^i)=(m+i)!x^i,\qquad i\geq0.
\end{equation}
If a graph $H$ contains $K_m$, then $a_k(H)=0$ for $k<m$.
Hence,
\begin{align}
 x^m\mathcal L_m\left(\frac{\sigma(H;x)}{x^m}\right)
 &=x^m\mathcal L_m\left(
      \sum_{i\geq0}a_{m+i}(H)x^i\right)\notag\\
 &=\sum_{i\geq0}(m+i)!a_{m+i}(H)x^{m+i}\notag\\
 &=\bar{\sigma}(H;x).
 \label{E4.8}
\end{align}
For $f(x),g(x)\in\RR[x]$, define the binary operation
$\diamond_m$ by
\begin{equation}\label{E4.9}
 f(x) \diamond_m g(x)
 =\sum_{j\geq0}\frac{[x(1+x)]^j}{j!(m+j)!}
                  (D^jf(x))(D^jg(x)).
\end{equation}
Similarly to the product $*$, the sum is finite. The following lemma relates this operation to the product $\ast$.
\begin{lemma}\label{L4.3}
Let $m\geq0$ be an integer and let $u(x),v(x)\in\RR[x]$.
If $f(x)=\mathcal L_m u(x)$ and $g(x)=\mathcal L_m v(x)$, then
\begin{equation}\label{E4.10}
 \mathcal L_m(u(x)\ast v(x))=f(x)\diamond_m g(x).
\end{equation}
\end{lemma}
\begin{proof}
Let $a,b$ be formal variables. We extend the operators
coefficientwise to formal power series in $a$ and $b$ with coefficients
in \(\RR[x]\). By equality \eqref{E4.5}, we have
\begin{align*}
 e^{ax}\ast e^{bx} =\sum_{j\geq0}\frac{x^j}{j!}
       a^je^{ax}b^je^{bx}
  =e^{(a+b)x}\sum_{j\geq0}\frac{(abx)^j}{j!}
 =e^{(a+b+ab)x}.
\end{align*}
By equality~\eqref{E4.7} and the binomial series, we obtain
\begin{align*}
 \mathcal L_m(e^{ax})
 =\sum_{i\geq0}\frac{(m+i)!}{i!}(ax)^i
 =m!\sum_{i\geq0}\binom{m+i}{i}(ax)^i
 =\frac{m!}{(1-ax)^{m+1}}.
\end{align*}
Hence
\begin{align*}
 \frac{\dd^j}{\dd x^j}\mathcal L_m(e^{ax})
 =\frac{(m+j)!a^j}{(1-ax)^{m+j+1}}.
\end{align*}

Let $A=x(1+x)$ and $C=(1-ax)(1-bx)$.
By equality~\eqref{E4.9}, we obtain
\begin{align*}
 (\mathcal L_m e^{ax})\diamond_m(\mathcal L_m e^{bx})
 &=\sum_{j\geq0}
   \frac{A^j}{j!(m+j)!}
   \frac{(m+j)!a^j}{(1-ax)^{m+j+1}}
   \frac{(m+j)!b^j}{(1-bx)^{m+j+1}}\\
 &=\sum_{j\geq0}\frac{(m+j)!}{j!}
             \frac{(abA)^j}{C^{m+j+1}}\\
 &=\frac{m!}{C^{m+1}}
   \sum_{j\geq0}\binom{m+j}{j}
                    \left(\frac{abA}{C}\right)^j\\
 &=\frac{m!}{C^{m+1}}
             \left(1-\frac{abA}{C}\right)^{-m-1}\\
 &=\frac{m!}{(C-abA)^{m+1}}.
\end{align*}
Since $C-abA =1-(a+b+ab)x$,
it follows that
\begin{align}\label{E4.11}
 (\mathcal L_m e^{ax})\diamond_m(\mathcal L_m e^{bx})
 =\frac{m!}{[1-(a+b+ab)x]^{m+1}}
 =\mathcal L_m(e^{(a+b+ab)x})
 =\mathcal L_m(e^{ax}\ast e^{bx}).
\end{align}
For each coefficient of $a^ib^\ell$, only the terms with
$j\leq\min(i,\ell)$ contribute to the differential sums.

Since $ \mathcal L_m(e^{ax}) =\sum_{i\geq0}\frac{a^i}{i!}\mathcal L_m(x^i)$ and 
\begin{align*}
 e^{ax}\ast e^{bx} =\bigg(\sum_{i\geq0}\frac{a^i}{i!}x^i\bigg)
   \ast
   \bigg(\sum_{\ell\geq0}\frac{b^\ell}{\ell!}x^\ell\bigg)
   =\sum_{i,\ell\geq0}\frac{a^ib^\ell}{i!\ell!}
   \,(x^i\ast x^\ell),
\end{align*}
the coefficient of \(a^ib^\ell\) on the right-hand side of equality
\eqref{E4.11} is
\begin{equation}\label{E4.12}
 \frac{1}{i!\ell!}\,
 \mathcal L_m(x^i\ast x^\ell).
\end{equation}
Similarly, the coefficient of $a^ib^\ell$ on the left-hand side
of equality \eqref{E4.11} is
\begin{equation}\label{E4.13}
 \frac{1}{i!\ell!}\,
 \bigl(\mathcal L_m(x^i)\diamond_m\mathcal L_m(x^\ell)\bigr).
\end{equation}
By equality \eqref{E4.11}, we have
\begin{align*}
 \frac{1}{i!\ell!}\,
 \mathcal L_m(x^i\ast x^\ell)
 =
 \frac{1}{i!\ell!}\,
 \bigl(\mathcal L_m(x^i)\diamond_m\mathcal L_m(x^\ell)\bigr).
\end{align*}
Multiplying both sides by $i!\ell!$ yields
\begin{align*}
 \mathcal L_m(x^i\ast x^\ell)
 =
 \bigl(\mathcal L_m(x^i)\bigr)\diamond_m
 \bigl(\mathcal L_m(x^\ell)\bigr),
\end{align*}
which is equality \eqref{E4.10} for $u=x^i$ and $v=x^\ell$.
The general case follows by bilinearity. This completes the proof.

\end{proof}
 
\section{Proof of the complete cutset theorem}\label{S5}
\hspace*{\parindent}
In this section, we first recall two stability results and construct
an auxiliary stable polynomial. Then we show that $\diamond_m$
preserves the property that all zeros belong to $[-1,0]$.

Let $\HH=\{z\in\mathbb C:\Impart z>0\}$.
A nonzero polynomial is called \emph{stable} if it does not vanish
when all its variables belong to $\HH$.
A nonzero univariate real polynomial is stable if and only if it
has only real zeros.
Let $e_k(t_1,\ldots,t_d)$ denote the elementary symmetric
polynomial of degree $k$ and $e_0=1$.

 The polarization result is obtained by polarizing only the variable $t$ in \cite[Proposition~2.4]{BB2009}.
\begin{proposition}[{\cite{BB2009}}]\label{P5.1}
If $h(z,t)=\sum_{k=0}^d h_k(z)t^k$ is stable, then
\begin{equation}\label{E5.1}
 \sum_{k=0}^d\frac{h_k(z)}{\binom dk}e_k(t_1,\ldots,t_d)
\end{equation}
is stable, where $z$ may denote one variable or a list of variables.
\end{proposition}

Another essential result is the univariate case of the sufficiency
direction of \cite[Theorem~1.3]{BB2010}.
\begin{theorem}[{\cite{BB2010}}]\label{t5.1}
Let $T=\sum_{j=0}^dQ_j(x)D^j$, where $Q_j\in\RR[x]$ and
$D=\dd/\dd x$. If its signed symbol
\begin{equation}\label{E5.2}
 S(x,z)=\sum_{j=0}^dQ_j(x)(-z)^j
\end{equation}
is stable, then $Tg$ is zero or has only real zeros whenever
$g$ has only real zeros.
\end{theorem}

The following lemma provides the auxiliary polynomial required for our proof.
\begin{lemma}\label{L5.1}
For integers $d,m\geq0$, the polynomial
\begin{equation}\label{E5.3}
 H_{d,m}(z;t_1,\ldots,t_d)
 =\sum_{j=0}^d\frac{z^j}{(m+j)!}e_{d-j}(t_1,\ldots,t_d)
\end{equation}
is stable.
\end{lemma}
\begin{proof}
For $d=0$, the polynomial is the nonzero constant $1/m!$.
For $d\geq1$,
by the explicit formula for the Laguerre polynomial
\cite[Eq.~(18.5.12)]{OL2010}
\begin{align*}
 L_d^{(m)}(s)
 =\sum_{j=0}^d\binom{d+m}{d-j}\frac{(-s)^j}{j!},
\end{align*}
we obtain
\begin{align}
 H_{d,m}(z;t,\ldots,t)
 =\sum_{j=0}^d\binom dj\frac{z^jt^{d-j}}{(m+j)!}
 =\frac{d!}{(m+d)!}\,t^dL_d^{(m)}(-z/t),
 \label{E5.4}
\end{align}
where the right-hand side is a polynomial in \(z\) and \(t\).

For $m\geq0$, the zeros $\rho_1,\rho_2,\ldots,\rho_d$ of $L_d^{(m)}(s)$
are positive \cite[Eq.~(18.16.9)]{OL2010}.
Since its leading coefficient is $(-1)^d/d!$, we have
\begin{align}\label{w1}
 L_d^{(m)}(s)
 =\frac{(-1)^d}{d!}\prod_{i=1}^d(s-\rho_i).
\end{align}
Substituting \eqref{w1} into~\eqref{E5.4}, we obtain
\begin{align*}
 H_{d,m}(z;t,\ldots,t)
 =\frac1{(m+d)!}\prod_{i=1}^d(z+\rho_i t).
\end{align*}
For $z,t\in\HH$, each factor satisfies
\begin{align*}
 \Impart(z+\rho_i t)
 =\Impart z+\rho_i\Impart t>0.
\end{align*}
Thus the diagonal polynomial is stable.

Its polarization in the variable $t$ is
\begin{align*}
 \sum_{j=0}^d
 \frac{\binom dj}{\binom d{d-j}}
 \frac{z^j}{(m+j)!}e_{d-j}(t_1,\ldots,t_d)
 =H_{d,m}(z;t_1,\ldots,t_d).
\end{align*}
By Proposition \ref{P5.1}, $H_{d,m}$ is stable.
This completes the proof.
\end{proof}

We next prove the stability of the signed symbol associated with
$\diamond_m$ using Lemma~\ref{L5.1}.
\begin{theorem}\label{T5.2}
Let $m\geq0$ be an integer. If $f(x)$ and $g(x)$ are nonzero real
polynomials whose zeros belong to $[-1,0]$, then $f(x) \diamond_m g(x)$ is nonzero and all its zeros belong to $[-1,0]$.
\end{theorem}
\begin{proof}
By bilinearity, we may assume that $f(x)$ and $g(x)$ are monic.
If either is constant, the result follows from~\eqref{E4.9}.
Assume that $d=\deg f(x)\geq1$ and $e=\deg g(x)\geq1$ and let
\begin{align}\label{l1}
 f(x)=\prod_{i=1}^d(x+a_i),\qquad 0\leq a_i\leq1.
\end{align}
Set $A(x)=x(1+x)$ and consider the differential operator
\begin{align*}
 T_f=\sum_{j=0}^d
       \frac{A(x)^jf^{(j)}(x)}{j!(m+j)!}D^j,
\end{align*}
where $f^{(j)}(x)$ is the $j$-th derivative of $f(x)$.
By equality \eqref{E4.9}, we have $T_f(g(x))=f(x) \diamond_m g(x)$ and its signed symbol is
\begin{equation}\label{E5.5}
 S_f(x,z)=\sum_{j=0}^d
          \frac{[-A(x)z]^j}{j!(m+j)!}f^{(j)}(x).
\end{equation}

For $x,z\in\HH$, let
\begin{equation}\label{E5.6}
 t_i=-\frac{x+a_i}{x(1+x)}
     =-\frac{a_i}{x}-\frac{1-a_i}{x+1},
 \qquad 1\leq i\leq d.
\end{equation}
Both $-1/x$ and $-1/(x+1)$ have positive imaginary part.
The coefficients $a_i$ and $1-a_i$ are nonnegative and sum to one,
so $t_i\in\HH$ for $a_i\in[0,1]$.

Differentiating equality \eqref{l1}, we obtain
\begin{align}\label{l2}
 \frac{f^{(j)}(x)}{j!}
 &=e_{d-j}(x+a_1,\ldots,x+a_d) \notag\\
 &=[-A(x)]^{d-j}e_{d-j}(t_1,\ldots,t_d).
\end{align}
 
Substituting \eqref{l2} into~\eqref{E5.5}, we obtain
\begin{align*}
 S_f(x,z)
 &=\sum_{j=0}^d
   \frac{[-A(x)z]^j[-A(x)]^{d-j}}{(m+j)!}
       e_{d-j}(t_1,\ldots,t_d)\\
 &=[-A(x)]^dH_{d,m}(z;t_1,\ldots,t_d).
\end{align*}
Since $A(x)\neq0$ for $x\in\HH$, Lemma~\ref{L5.1}
gives $S_f(x,z)\neq0$ for $x,z\in\HH$.
Thus $S_f$ is stable.
 
The polynomial $f(x) \diamond_m g(x)$ is zero or has only real zeros by Theorem \ref{t5.1}. It remains to show that this polynomial is nonzero and has no
zero outside $[-1,0]$.
For $x>0$, all nonzero derivatives of $f(x)$ and $g(x)$ are positive and $A(x)>0$. Hence each nonzero summand in equality~\eqref{E4.9}
is positive. In particular, the summand with $j=0$ is
$f(x)g(x)/m!>0$. Therefore $f(x)\diamond_m g(x)$ is nonzero and
has no positive zero.

For $x<-1$, let
\begin{align*}
 g(x)=\prod_{i=1}^e(x+b_i),\qquad 0\leq b_i\leq1.
\end{align*}
All factors $x+a_i$ and $x+b_i$ are negative.
Since $f^{(j)}(x)/j!$ is a sum of products
of $d-j$ such factors and similarly for $g^{(j)}(x)/j!$, we obtain
\begin{align*}
 \sgn f^{(j)}(x)=(-1)^{d-j},\qquad
 \sgn g^{(j)}(x)=(-1)^{e-j},
 \qquad 0\leq j\leq\min(d,e).
\end{align*}
Clearly $A(x)>0$, so each nonzero summand in equality~\eqref{E4.9} has sign
\begin{align*}
 (-1)^{d-j}(-1)^{e-j}=(-1)^{d+e}.
\end{align*}
Hence, these summands cannot cancel and the summand with $j=0$
is nonzero. Thus $f(x)\diamond_m g(x)$ has no zero in $(-\infty,-1)$.
It follows that all its zeros belong to $[-1,0]$.
This completes the proof.
\end{proof}

We next prove one of the main results of this paper.

\begin{proof}[Proof of Theorem~\ref{T1.2}]
If either graph has no vertices, the result is immediate.
Assume that both graphs are nonempty.

Since $G\cap H=K_m$, both $G$ and $H$ contain $K_m$.
Thus $\bar{\sigma}(G;x)$ and $\bar{\sigma}(H;x)$ are divisible
by $x^m$. By Lemma~\ref{L4.1}, the polynomials
\begin{align*}
 f(x)&=\frac{\bar{\sigma}(G;x)}{x^m}
      =\mathcal L_m\left(\frac{\sigma(G;x)}{x^m}\right),\\
 g(x)&=\frac{\bar{\sigma}(H;x)}{x^m}
      =\mathcal L_m\left(\frac{\sigma(H;x)}{x^m}\right)
\end{align*}
are nonzero and all their zeros belong to $[-1,0]$,
where the equalities follow from equality~\eqref{E4.8}.

The graph $G\cup H$ also contains $K_m$.
By equality~\eqref{E4.8}, Theorem~\ref{T4.2} and Lemma~\ref{L4.3},
we obtain
\begin{align*}
 \bar{\sigma}(G\cup H;x)
 &=x^m\mathcal L_m\left(
       \frac{\sigma(G\cup H;x)}{x^m}\right)\\
 &=x^m\mathcal L_m\left(
       \frac{\sigma(G;x)}{x^m}
       \ast
       \frac{\sigma(H;x)}{x^m}\right)\\
 &=x^m(f(x)\diamond_m g(x)).
\end{align*}
By Theorem~\ref{T5.2}, $f(x)\diamond_m g(x)$ is nonzero and all its
zeros belong to $[-1,0]$.
Multiplication by $x^m$ adds only zeros at zero.
Hence all zeros of $\bar{\sigma}(G\cup H;x)$ belong to $[-1,0]$.
By Lemma~\ref{L4.1}, $w(G\cup H;x)\in\RZ$.
This completes the proof.
\end{proof}

%\section*{References}

\end{document}